\documentclass[amssymb,amsfonts]{amsart}
\usepackage{amssymb, latexsym}
\def\be#1{\begin{equation} \label{#1}}
\def\bi{\begin{itemize}}
\def\bs{\begin{split}}
\def\es{\end{split}}
\def\rr{{\mathbb R}}
\def\ba{\begin{align}}
\def\bas{\begin{align*}}
\def\philambda{\phi^{(\lambda)}}
\def\phihigh{\phi^{\text{high}}}
\def\philow{\phi^{\text{low}}}
\def\rtwo{{\bf R}^2}
\def\ea{\end{align}}
\def\eas{\end{align*}}

\newcommand{\rn}{{\mathbb R}^n}
\def\Re{{\hbox{Re}}}
\def\rr{{\bf R}}
\def\fatD{{\langle \nabla \rangle}}

\def\Xdoohp{{X^{\delta}_{1, \frac{1}{2} +}}}
\def\Xoohp{X_{1, \frac{1}{2}+ } }
\def\Xzohp{X_{0, \frac{1}{2}+}  }
\def\Xdzohp{\Xzohp^{\delta}}

\def\termone{{\text{Term}_{{1}} }}
\def\termtwo{{\text{Term}_{{2}} }}

\def\textbf#1{{\bf #1}}

\theoremstyle{plain}
  \newtheorem{theorem}[subsection]{Theorem}
  
  \newtheorem{proposition}[subsection]{Proposition}
  \newtheorem{lemma}[subsection]{Lemma}

\theoremstyle{remark}

\numberwithin{equation}{section}
\theoremstyle{definition}

\include{psfig}

\begin{document}

\title[A.C. Laws and Global Rough Solutions for NLS]{Almost Conservation Laws \\
and  Global Rough Solutions\\
 to a Nonlinear Schr\"odinger Equation}

\vspace{-0.3in}

\author{J. Colliander}
\thanks{J.E.C. was supported in part by N.S.F. Grant DMS 0100595.}
\address{\small University of Toronto}

\author{M. Keel}
\thanks{M.K. was supported in part by N.S.F. Grant DMS 9801558.}
\address{\small University of Minnesota, Minneapolis}

\author{G. Staffilani}
\thanks{G.S. was supported in part by N.S.F. Grant 
DMS 0100375 and
grants from Hewlett and Packard and the Sloan Foundation.}
\address{\small Brown University and Stanford University}

\author{H. Takaoka}
\address{\small Hokkaido University}
\thanks{H.T. was supported in part by J.S.P.S. Grant No. 13740087}

\author{T. Tao}
\thanks{T.T. is a Clay Prize Fellow and was supported in part by
a grant from the Packard Foundation.}
\address{\small University of California, Los Angeles}


\begin{abstract}
We prove an ``almost conservation law"
 to obtain global-in-time well-posedness for the cubic, defocussing 
nonlinear Schr\"odinger equation in $H^s(\rr^n)$ when $n = 2,3$ and 
$s > \frac{4}{7}, \frac{5}{6}$, respectively.
\end{abstract}

\maketitle

\section{Introduction and Statement of Results}

We study the following initial value problem for a defocussing nonlinear
Schr\"odinger equation, 
\begin{align}
\label{nls}
i \partial_t \phi(x,t) + \Delta \phi (x,t)& = |\phi(x,t)|^2 \phi(x,t) \quad
x \in \rr^n, t \geq 0   \\
\phi(x,0) & = \phi_0(x) \; \in  H^s(\rr^n) \label{nlsdata}  
\end{align}
when $n = 2,3$.  Here $H^s(\rr^n)$ denotes the usual inhomogeneous 
Sobolev space.  Our goal is to loosen the regularity requirements on the
initial data which ensure global-in-time solutions.  In particular, 
we aim to extend the global theory to certain infinite energy initial
data.

It is known 
\cite{cw:local} 
that \eqref{nls}-\eqref{nlsdata} is well-posed locally in time 
when $n=2,3$ and $s >0,
\frac{1}{2}$ respectively\footnote{In addition, there are local
in time solutions from $L^2, H^{\frac{1}{2}}$ data 
when $n=2,3$, respectively.  However, it is not yet known whether 
the time interval of existence for such solutions depends only 
on the data's Sobolev norm.  For example, the $L^2$ conservation
law \eqref{l2conservation} does not yield the widely conjectured
result of global in time
solutions on $\rr^{2+1}$ from $L^2$ initial data.}.  
In addition,  these local solutions
enjoy $L^2$ conservation;
\begin{align}
\label{l2conservation}
||\phi(\cdot, t)||_{L^2(\rn)} & = ||\phi_0(\cdot)||_{L^2(\rn)}
\end{align}
and the $H^1(\rr^n)$ solutions have the following conserved energy,
\begin{align}
\label{energy}
E(\phi)(t) & \equiv \int_{\rn} \frac{1}{2} 
|\nabla_ x\phi(x,t)|^2 + \frac{1}{4} |\phi(x,t)|^4\ dx \quad = \;
E(\phi)(0).
 \end{align}

Together,  energy conservation and the local-in-time theory immediately
yield global-in-time well-posedness of \eqref{nls}-\eqref{nlsdata} 
from data in $H^s(\rn)$ when
$s \geq 1$, and $n = 2,3$. It is conjectured that 
\eqref{nls}-\eqref{nlsdata} is in fact  globally well-posed in time from all
data included in the local theory.  
The obvious impediment to claiming global-in-time solutions in $H^s$, 
with $0< s <1$, is the lack of any applicable conservation law.

The first argument extending the lifespan of 
rough solutions to \eqref{nls}-\eqref{nlsdata} 
in a range $s_0 < s < 1$
was given in \cite{bourg1} (see also \cite{bourgbook}).  In 
what might be called a  
``Fourier truncation" approach, Bourgain observed that from
the point of view of regularity, the high frequency component of 
the solution $\phi$ is well-approximated by the corresponding 
linear evolution of the data's high frequency component.  More
specifically:  one makes a first approximation to the solution
for a small time step by evolving the high modes 
linearly, and the low modes according to the nonlinear
flow for which one has energy conservation.  The correction term
one must add to match this approximation with the actual solution
is shown to have finite energy.  This correction is added to
the low modes as data for the nonlinear evolution during
the next time step, where the high modes are again evolved linearly.  
For $s > \frac{3}{5}$, one 
can repeat this procedure to an arbitrarily large time provided the
distinction between ``high" and ``low" frequencies is made at
sufficiently large frequencies.

The argument in \cite{bourg1} has been applied to other 
subcritical initial value problems with sufficient smoothing in their
principal parts.  (See \cite{bourgbook}, \cite{cst:kdv},
\cite{FLP}, \cite{kpv:wave}, \cite{takaoka:dnls}, 
and \cite{hideo:kp}).
It is important to note that the Fourier truncation method demonstrates more
than  just rough data global existence.  Indeed,  write  $S_t^{NL}$  for the
nonlinear flow\footnote{That is,   $S_t^{NL}(\phi_0)(x) \, = \phi(x,t)$, 
where
$\phi, \phi_0$ as in   \eqref{nls}-\eqref{nlsdata}.} of
\eqref{nls}-\eqref{nlsdata}, and  let $S_t^{L}$ denote the corresponding
{\em{linear}} flow.  The Fourier truncation method shows then that for $s > \frac{3}{5}$ 
and for  all $t \in [0, \infty)$, 
\begin{align} 
\label{smoothing} 
S_t^{NL}\phi_0 - S^L_t \phi_0 & \in H^{1}(\rr^2). 
\end{align} 
Besides being part of the conclusion, the smoothing 
property \eqref{smoothing} seems to be a crucial constituent of
the Fourier truncation argument itself. 

In this paper we will use a modification of the above arguments,
originally put forward to analyze equations where the smoothing property 
\eqref{smoothing} is not available because it is either false
(e.g. Wave maps \cite{wavemaps}\footnote{See the appendix of \cite{mkg}
for the failure of \eqref{smoothing} for Wave Maps.}) or 
simply not known (e.g. Maxwell-Klein-Gordon equations \cite{mkg}, for
which we suspect \eqref{smoothing} is false). In this
``almost conservation law" approach,  one  controls the growth in time of
a rough solution by monitoring the energy of a certain smoothed out version
of the solution.  It can be shown that the energy of the smoothed solution
is ``almost conserved" as time passes, and controls the solution's 
sub-energy Sobolev norm.  In proving the almost conservation law
for the i.v.p. \eqref{nls}-\eqref{nlsdata}, we shall use only the
linear estimates presented in  \cite{bourg1}, \cite{bourgbook}.  Implicitly, we
also use the
view of  \cite{bourg1} that the energy at high frequencies does not move 
rapidly to low frequencies.

The almost conservation approach to global rough solutions has 
proven to be quite robust \cite{wavemaps}, \cite{mkg}, \cite{firstkdv}, \cite{iteamIV}, and has  
been improved
significantly by  adding additional 
``correction" terms to the original
almost conserved energy functional.  As a result, one obtains even 
stronger
bounds on the growth of the solution's rough norm, and at least 
in some cases
sharp global well-posedness results \cite{iteamderiv}, \cite{iteamI}, 
\cite{iteamIII}.   

The aims of this paper are three-fold:   first and most 
obviously, an improved understanding of the evolution properties of 
rough solutions 
of \eqref{nls}-\eqref{nlsdata};  second, 
the almost conservation law approach is presented in a relatively 
straightforward setting; 
and third,  
we can directly compare this
almost conservation law approach to the Fourier
cut-off technique, since both approaches apply to the semilinear
Schr\"odinger initial value problem.  Our main result is the following:

\begin{theorem}  \label{maintheorem}
The initial value problem \eqref{nls}-\eqref{nlsdata} is globally-well-posed from data 
$\phi_0 \in H^s(\rr^n), \; n = 2,3$ when $s > \frac{4}{7}, \frac{5}{6}$ 
respectively.  
\end{theorem}

By ``globally-well-posed", we mean that given 
data $\phi_0 \in H^s(\rr^n)$ as above, and any time $T > 0$,  there
is a unique solution to \eqref{nls}-\eqref{nlsdata}
\begin{align}
\phi(x,t) \in C([0,T]; H^s(\rr^n)) \label{summary}
\end{align}
which depends continuously in \eqref{summary} upon $\phi_0 \in H^s(\rr^n)$. The 
polynomial bounds we obtain for the growth of $||\phi||_{H^s(\rr^n)}(t)$  
are contained in \eqref{boundhsnorm}, \eqref{donepolybound}, and \eqref{polynomial3d} 
below.

Theorem \ref{maintheorem} extends to some extent the work in
\cite{bourg1, bourgbook}  where global well-posedness was shown when 
$s >\frac{3}{5}, \frac{11}{13}$ and $n=2,3$ respectively. In a different
sense, the result here is weaker than the results of \cite{bourg1, bourgbook}
as we obtain no information whatsoever along the lines of
\eqref{smoothing}. 

In a later paper, we hope to extend Theorem \ref{maintheorem} to still 
rougher data, using the additional cancellation terms mentioned above, and
the multilinear estimates contained in \cite{fancynls}. 

In Section \ref{section:notation} below we present some notation
and linear estimates that are used in our proofs.  Sections 
\ref{section:2d}, \ref{section:3d} present the almost conservation 
laws and proofs of Theorem \ref{maintheorem} in space dimensions
two and three, respectively.

\section{Estimates, Norms, and Notation}
\label{section:notation}
Given  $A,B \geq 0$, we write
$A \lesssim B$ to mean that for some universal constant $K > 2$,
$A \leq K\cdot B$.  We write $A \thicksim B$ when both $A \lesssim B$ 
and $B \lesssim A$.  The notation $A \ll B$ denotes $B > K\cdot A$. 

We write $\langle A \rangle \equiv (1 + A^2)^{\frac{1}{2}}$, and 
$\fatD$ for the operator with Fourier multiplier $(1 + |\xi|^2)^{\frac{1}{2}}$.
The symbol $\nabla$ will denote the spatial gradient.

We will use the weighted Sobolev norms,
(see
\cite{rauch-reed, beals, bourgainnorm, kman_machedon}), 
\begin{align}
\label{norms}
||\psi||_{X_{s, b}}  & \equiv ||\langle\xi\rangle^s 
\langle \tau -
|\xi|^2 \rangle^{b} \tilde{\psi}(\xi, \tau) ||_{L^2(\rr^n \times \rr)}.
\end{align}  
Here $\tilde{\psi}$ is the space-time Fourier transform of $\psi$.
We will need local-in-time estimates in terms of truncated versions
of the norms \eqref{norms},
\begin{align}
\label{truncatednorms}
||f||_{X_{s,b}^{\delta}} & \equiv \inf_{\psi = f \text{on}
[0, \delta]} ||\psi||_{X_{s, b}^{\delta}}.
\end{align}
We will often use the notation $\frac{1}{2}+ \equiv \frac{1}{2} + \epsilon$ for some
universal $0 < \epsilon \ll 1$. Similarly, we shall write
$\frac{1}{2}- \equiv \frac{1}{2} - \epsilon$, and 
$\frac{1}{2} -- \equiv \frac{1}{2} - 2 \epsilon$.  

Given Lebesgue space exponents $q,r$ and a function $F(x,t)$ on 
$\rr^{n+1}$, we write
\begin{align}
\label{mixedlebesgue}
||F||_{L^q_tL^r_x(\rr^{n+1})} & \equiv \left( \int_{\rr}
\left( \int_{\rr^n} |F(x,t)|^r dx \right)^{\frac{q}{r}} dt
\right)^{\frac{1}{q}}.
\end{align}
This norm will be shortened to $L^q_tL^r_x$ for readability, 
or to $L^r_{x,t}$ when $q=r$.  

We will need Strichartz-type estimates \cite{yajima, gv, endpoint} involving the
spaces \eqref{mixedlebesgue}, \eqref{norms}.  We will call a 
pair of exponents $(q,r)$  {\em{ Schr\"odinger admissible}} for
$\rr^{n+1}$ when $q,r \geq 2$, $(n,q) \neq (2,2)$, and 
\begin{align}
\label{sa}
\frac{1}{q} + \frac{n}{2r} & = \frac{n}{4}.
\end{align}  
For a Schr\"odinger admissible pair $(q,r)$ we have
what we will call the $L^q_tL^r_x$ Strichartz estimate,
\begin{align}
\label{strichartz}
||\phi||_{L^q_t L^r_x(\rr^{n+1})} & \lesssim ||\phi||_{\Xzohp}.
\end{align}

Finally, we will need a refined version of these estimates
due to Bourgain \cite{bourg1}.

\begin{lemma}  \label{bourgainstrichartz}  
Let $\psi_1, \psi_2 \in \Xdzohp$ be  
supported on spatial frequencies $|\xi| \thicksim N_1, N_2$, respectively.
Then for   $N_1 \leq N_2$, one has 
\begin{align} 
\label{bstrich} 
||\psi_1 \cdot \psi_2||_{L^2([0, \delta] \times \rr^2)}  & 
\lesssim \left(  \frac{N_1}{N_2}
\right)^{\frac{1}{2}} ||\psi_1||_{\Xdzohp}   ||\psi_2||_{\Xdzohp}.  
\end{align} 
In addition, \eqref{bstrich} holds (with the same proof) if we replace the 
product $\psi_1 \cdot \psi_2$ on 
the left with either $\overline{\psi}_1 \cdot \psi_2$
or $\psi_1 \cdot \overline{\psi}_2$. 
\end{lemma}

\section{Almost conservation and Proof of Theorem 
\ref{maintheorem} in $\rr^2$}\label{section:2d}

For rough initial data, \eqref{nlsdata} with $s < 1$, the energy
is infinite, and so the conservation law \eqref{energy}
is meaningless.  
Instead, Theorem \ref{maintheorem} rests on the fact that  a smoothed 
version of the solution \eqref{nls}-\eqref{nlsdata} has a finite 
energy which 
is {\em almost} conserved in time.  We express this `smoothed version' as 
follows.

Given $s < 1$ and a parameter $N \gg 1$, define the multiplier operator
\begin{align}
\label{Ioperator}
\widehat{I_N f}(\xi) & \equiv m_N(\xi) \hat{f}(\xi), 
\end{align}
where the multiplier $m_N(\xi)$ is smooth, radially symmetric, nonincreasing
in $|\xi|$ and 
\begin{align}
\label{Iproperties}
m_N(\xi) & =  \begin{cases}
1 & |\xi| \leq N \\
\left( \frac{N}{|\xi|} \right)^{1 - s} & |\xi| \geq 2N.
\end{cases}
\end{align}
For simplicity, we will eventually drop the $N$ from the notation,
writing $ I$ and $m$ for \eqref{Ioperator} and 
\eqref{Iproperties}.  
Note that for solution and initial data $\phi, \, \phi_0$ of 
\eqref{nls}, \eqref{nlsdata},  the quantities 
$||\phi||_{H^s(\rn)}(t)$ and $E(I_N \phi)(t)$ (see \eqref{energy})
can be compared,
\begin{align} \label{boundenergy}
E(I_N \phi)(t) & \leq \left(N^{1 - s} 
||\phi(\cdot, t)||_{\dot H^s(\rn)} \right)^2 + 
||\phi(t, \cdot)||_{L^4(\rn)}^4, \\
\label{boundhsnorm}
||\phi(\cdot, t)||_{H^s(\rn)}^2 & \lesssim E(I_N \phi)(t) + 
||\phi_0||_{L^2(\rn)}^2.
\end{align}
Indeed, the $\dot H^1(\rn)$ component of the left hand side of
\eqref{boundenergy} is bounded by the right side by using the 
definition of $I_N$ and
by considering separately those frequencies
$|\xi| \leq N$ and $|\xi| \geq N$. The $L^4$ component of 
the energy in \eqref{boundenergy} is bounded
by the right hand side of \eqref{boundenergy} by using (for example) the 
H\"ormander-Mikhlin multiplier theorem.  The bound \eqref{boundhsnorm}
follows quickly from \eqref{Iproperties} and 
$L^2$ conservation \eqref{l2conservation} by considering separately
the $\dot{H}^s(\rr^n)$ and $L^2(\rr^n)$ components of the left hand
side of \eqref{boundhsnorm}.

To prove Theorem \ref{maintheorem}, we may
assume that $\phi_0  \in C^{\infty}_0(\rr^n)$, and show
that the resulting global-in-time solution
grows at most polynomially
in the $H^s$ norm,
\begin{align}
\label{polynomialgrowth}
||\phi(\cdot, t)||_{H^s(\rn)} & \leq C_1 t^M + C_2, 
\end{align}
where the constants $C_1, C_2, M$ depend only 
on $||\phi_0||_{H^s(\rr^n)}$ and
not on higher regularity norms of the smooth data.
Theorem \ref{maintheorem} follows immediately
from \eqref{polynomialgrowth}, the 
local-in-time theory \cite{cw:local}, and a standard density argument. 

By \eqref{boundhsnorm}, it suffices to show 
\begin{align}
\label{polynomialgrowthII}
E(I_N\phi)(t) & \lesssim (1+t)^M.
\end{align}
for some $N = N(t)$. (See \eqref{here'sN}, \eqref{donepolybound}
below for the definition of $N$ and the growth rate $M$ we eventually
establish.)  The following proposition, which is one of  the two 
main estimates of this paper (see also Proposition \ref{almostconservation3d}),
represents an ``almost conservation law'' of 
the title and will yield \eqref{polynomialgrowthII} in space dimension $n=2$.

\begin{proposition}  \label{almostconservation}
Given $s > \frac{4}{7}, N \gg 1,$ and initial data 
$\phi_0 \in C^{\infty}_0(\rtwo)$ (see preceeding
remark) with $E(I_N \phi_0) \leq 1$, then there exists a $ \delta =
\delta(||\phi_0||_{L^2(\rtwo)}) > 0$  so that the solution
\begin{align*}
\phi(x,t) & \in C([0,\delta], H^s(\rtwo))
\end{align*}
of \eqref{nls}-\eqref{nlsdata} satisfies
\begin{align}
\label{increment}
E(I_N\phi)(t) & =  E(I_N\phi)(0) + O(N^{- \frac{3}{2}+}),
\end{align}
for all $t \in [0, \delta]$.
\end{proposition}
 
\noindent{\bf Remark}:  Equation \eqref{increment} asserts that 
$I_N\phi$, though not a solution of the 
nonlinear problem \eqref{nls}, enjoys something akin to energy conservation.
If one could replace the increment $N^{- \frac{3}{2}+}$ in $E(I_N\phi)$ on the right side of 
\eqref{increment} with $N^{- \alpha}$ for some $\alpha > 0$, one could repeat the 
argument we give
below to prove global well posedness of 
\eqref{nls}-\eqref{nlsdata} for all $s > \frac{2}{2+\alpha}$.  
In particular, if $E(I_N \phi)(t)$ is conserved (i.e. $\alpha = \infty$),
one could show that \eqref{nls}-\eqref{nlsdata} is globally well-posed when $s > 0$.
 
We first show that Proposition \ref{almostconservation}
implies \eqref{polynomialgrowthII}.    Note that the initial value problem
here has a scaling symmetry, and  is $H^s$-subcritical when 
$1 > s > 0, \frac{1}{2}$
and $n = 2,3$, respectively.   That
is, if $\phi$ is a solution to \eqref{nls}, so too
\begin{align}
\label{scaling}
\philambda(x,t) & \equiv \frac{1}{\lambda} 
\phi(\frac{x}{\lambda}, \frac{t}{\lambda^2}).
\end{align}
Using \eqref{boundenergy}, the following energy
can be made
arbitrarily small by taking $\lambda$ large,
\begin{align}
E(I_N \philambda_0) & \leq \left( (N^{2 - 2s} ) \lambda^{-2s} + \lambda^{-2}
\right) \cdot (1 + ||\phi_0||_{H^s(\rr^2)})^4 \\
& \leq C_0(N^{2 - 2s} \lambda^{-2s} ) \cdot  (1 + ||\phi_0||_{H^s(\rr^2)})^4. 
\label{makeitsmall}
\end{align}
Assuming $N \gg 1$ is given\footnote{The parameter $N$ will be chosen shortly.}, we
choose our scaling parameter $\lambda = \lambda(N, ||\phi||_{H^s(\rr^2)})$
\begin{align}
\label{chooselambda}
\lambda & = N^{\frac{1 -s}{s}} \left( \frac{1}{2C_0} \right)^{- \frac{1}{2s}}
\cdot \left(1 + ||\phi_0||_{H^s(\rr^2)}\right)^{\frac{2}{s}}
\end{align}
so that $E(I_N\philambda_0) \leq \frac{1}{2}$. 
We may now apply Proposition \ref{almostconservation} to the scaled initial
data $\philambda_0$, and in fact may reapply this Proposition until the size
of $E(I_N\philambda)(t)$ reaches $1$, that is 
at least  $C_1 \cdot N^{\frac{3}{2}-}$ times. Hence 
\begin{align}
\label{scaledbound}
E(I_N \philambda)(C_1N^{\frac{3}{2}-}\delta) &  \thicksim 1.
\end{align}

Given any $T_0 \gg 1$, we establish the polynomial growth \eqref{polynomialgrowthII}
from 
\eqref{scaledbound}
by first choosing our parameter $N \gg 1$ so that 
\begin{equation}\label{here'sN}
T_0 \thicksim \frac{N^{\frac{3}{2}-}}{\lambda^2} 
C_1 \cdot \delta \thicksim N^{\frac{7s-4}{2s}-},
 \end{equation}
where we've kept in mind  \eqref{chooselambda}.  Note the exponent 
of $N$ on the right of \eqref{here'sN} is
positive provided $s > \frac{4}{7}$, hence the definition of $N$
makes sense for arbitrary $T_0$.  In two space dimensions, 
\begin{align*}
E(I_N \phi)(t) & = \lambda^2 E(I_N \philambda)(\lambda^2 t).
\end{align*}
We use \eqref{chooselambda}, \eqref{scaledbound}, and \eqref{here'sN} to 
conclude that for $T_0 \gg  1$, 
\begin{align}
\label{donepolybound}
E(I_N\phi)(T_0) & \leq C_2 T_0^{\frac{1-s}{\frac{7}{4}s - 1} +},
\end{align}
where $N$ is chosen as in \eqref{here'sN} and
$C_2 = C_2(||\phi_0||_{H^s(\rr^2)}, \delta)$.  Together with
\eqref{boundhsnorm}, the bound \eqref{donepolybound} establishes
the desired polynomial bound \eqref{polynomialgrowth}.

It remains then to prove Proposition \ref{almostconservation}.
We will need the following modified version of the usual local
existence theorem, wherein we control for small times  
the smoothed solution in the $X_{1, \frac{1}{2}+}^{\delta}$ norm.

\begin{proposition} \label{modifiedlocalexist} Assume 
$\frac{4}{7} < s < 1$ and we are given data
for the problem \eqref{nls}-\eqref{nlsdata} with $E(I\phi_0) \; \leq \; 1$. 
Then there is a constant $\delta =
\delta(||\phi_0||_{L^2(\rtwo)})$ so that the solution $\phi$ obeys
the following bound on the time interval $[0, \delta]$,
\begin{align}
\label{localbounds}  
||I\phi||_\Xdoohp 
& \lesssim 1.
\end{align}
\end{proposition}

\begin{proof} 
We mimic the typical iteration argument showing local existence.  
We will need the following three estimates involving the $X_{s,
\delta}$ spaces \eqref{norms} and 
functions $F(x,t), f(x)$.  (Throughout this section, the implicit
constants in the notation $\lesssim$ are independent of $\delta$.)    
\begin{align}
\|S(t) f\|_{\Xdoohp} & \lesssim \|f\|_{H^1(\rr^2)}, \label{one} \\ 
\left\| \int_0^t S(t - \tau) F(x, \tau) d \tau \right\|_{\Xoohp}
& \lesssim  \|F\|_{X^\delta_{1, -\frac{1}{2} +}}, \label{two} \\
\|F\|_{X^{\delta}_{1, -b}} & \lesssim  \delta^P 
\|F\|_{X^\delta_{1, -\beta}},
\label{three}
\end{align}
where in \eqref{three} we have $0 < \beta < b < \frac{1}{2}$, and 
$P = \frac{1}{2} (1 - \frac{\beta}{b}) > 0.$  The bounds \eqref{one}, \eqref{two}
are analogous to estimates (3.13), (3.15) in \cite{kpv:bilinear}.  As
for \eqref{three},
by duality it suffices to show
\begin{align*}
||F||_{X_{-1, \beta}^{\delta}} & \lesssim \delta^P
||F||_{X_{-1,b}^{\delta}}.
\end{align*}
Interpolation\footnote{The argument here actually involves 
Lemma 3.2 of \cite{kpv:bilinear}.  We thank S. Selberg for pointing
this out to us.} gives 
\begin{align*}
||F||_{X_{-1, \beta}^\delta}
& \lesssim ||F||_{X_{-1, 0}^\delta} ^{(1 -
\frac{\beta}{b})-} \cdot ||F||_{X_{-1, b}^\delta}
^{\frac{\beta}{b}}.
\end{align*}
As $b \in (0, \frac{1}{2})$, arguing exactly as on page 771
of \cite{cst:kdv},
\begin{align*}
||F||_{X^\delta_{-1, 0}} & \lesssim
\delta^{\frac{1}{2}} ||F ||_{X^{\delta}_{-1, b}},
\end{align*}
and \eqref{three} follows.

Duhamel's principle and \eqref{one}- \eqref{three} give us
\begin{align}
||I\phi||_{\Xdoohp} 
& = \left\| S(t)(I \phi_0) + 
\int_0^t S(t - \tau) I(\phi \bar \phi \phi)(\tau) d\tau 
\right\|_{\Xdoohp} \nonumber  \\
& \lesssim ||I\phi_0||_{H^1(\rtwo)}+ || 
I(\phi \overline \phi \phi) 
||_{X^{\delta}_{1, -\frac{1}{2} +}} \nonumber \\
& \lesssim ||I\phi_0||_{H^1(\rtwo)} + \delta^{\epsilon} || 
I(\phi \overline{\phi} \phi) ||_{X^\delta_{
1, - \frac{1}{2}++} }, \label{readyforlove} 
\end{align}
where $-\frac{1}{2} + +$ is a real number slightly larger 
than $-\frac{1}{2}+$ and $\epsilon > 0$.  
By the definition of the restricted norm \eqref{truncatednorms},
\begin{align}
\label{switcheroo}
||I\phi||_{\Xdoohp} & \lesssim ||I\phi_0||_{H^1(\rtwo)} + 
\delta^{\epsilon} ||
I(\psi \overline{\psi} \psi)||_{X_{1, - \frac{1}{2}+ +} },
\end{align}
where the function $\psi$ agrees with $\phi$ for $t \in [0, \delta]$, and 
\begin{align}
\label{switchreason}
||I\phi||_{\Xdoohp} & \sim \quad ||I \psi ||_{X_{1, \frac{1}{2} + }}.
\end{align}
We will show shortly that 
\begin{align}
|| I(\psi \overline
 \psi \psi) ||_{X_{1,  -\frac{1}{2} ++ }} & 
\lesssim \quad || I \psi ||^3_{\Xoohp}.
\label{wholetthedogsout}
\end{align}
Setting then 
$Q(\delta) \,  \equiv \,||I \phi(t) ||_{\Xdoohp}$,
the bounds
\eqref{readyforlove}, \eqref{switchreason} and 
\eqref{wholetthedogsout} yield
\begin{align}
\label{weletthedogsout}
Q(\delta) & \lesssim ||I\phi_0||_{H^1(\rr^2)} + \delta^{\epsilon}
(Q(\delta))^3.
\end{align}
Note 
\begin{align}
||I\phi_0||_{H^1(\rr^2)} &  \lesssim 
(E(I\phi_0))^{\frac{1}{2}} + ||\phi_0||_{L^2(\rr^2)} 
\lesssim 1 + ||\phi_0||_{L^2(\rr^2)}. \label{datapain}
\end{align}
As $Q$ is continuous in the variable $\delta$, a bootstrap argument 
yields \eqref{localbounds} from \eqref{weletthedogsout}, \eqref{datapain}.

It remains to show \eqref{wholetthedogsout}.  Using the interpolation lemma of 
\cite{iteamIII}, 
it suffices to show
\begin{align}
\label{wltdo}
||\psi \bar{\psi} \psi ||_{X_{s, - \frac{1}{2} ++}} & 
\lesssim || \psi ||^3_{X_{s, \frac{1}{2} +}},
\end{align} 
for all $\frac{4}{7} < s < 1$. 
By duality and a ``Leibniz" rule\footnote{By this, we mean the operator $\langle D \rangle^s$
can be distributed over the product by taking Fourier transform and using
$\langle \xi_1 + \ldots \xi_4\rangle^s \lesssim \langle\xi_1 \rangle^s + 
\ldots \langle\xi_4\rangle^s$.},  \eqref{wltdo} follows from  
\begin{align} \label{dualwltdo}
\left| \int_{\rr} \int_{\rr^2} (\langle 
\nabla \rangle^s \phi_1) \overline{\phi_2} \phi_3 \phi_4 dx dt \right| & 
\lesssim ||\phi_1||_{X_{s, \frac{1}{2} +}} 
\cdot ||\phi_2||_{X_{s, \frac{1}{2} +}} \cdot
||\phi_3||_{X_{s, \frac{1}{2} +}} ||\phi_4||_{X_{0, \frac{1}{2}--}}.
\end{align}
Note that since the factors in the integrand on the left here will 
be taken in absolute value, the relative placement of 
complex conjugates is irrelevant.
Use H\"older's inequality on the left side of \eqref{dualwltdo}, taking 
the factors in, respectively, $L^4_{x,t}, L^4_{x,t}, L^6_{x,t}$ and
$L^3_{x,t}$.  Using a
Strichartz inequality, 
\begin{align*} 
||\langle \nabla \rangle^s \phi_1||_{L^4_{x,t}(\rr^{2+1})} & \lesssim || \langle
 \nabla \rangle^s \phi_1||_{\Xzohp} \\
& = || \phi_1 ||_{X_{s, \frac{1}{2}+}},
\end{align*}
and 
\begin{align*}
||\phi_2||_{L^4_{x,t}(\rr^{2+1})} & \lesssim ||\phi_2||_{\Xzohp} \\
& \lesssim ||\phi_2||_ {X_{s, \frac{1}{2}+}}.
\end{align*}
The bound for the third factor uses Sobolev embedding and the 
$L^6_tL^3_x$ Strichartz estimate,
\begin{align*}
||\phi_3||_{L^6_tL^6_x(\rr^{2+1})} & \lesssim ||\langle \nabla \rangle^{\frac{1}{3}}
\phi_3||_{L^6_tL^3_x(\rr^{2+1})} \\ 
& \lesssim ||\langle \nabla \rangle^{\frac{1}{3}} \phi_3||_{X_{0, \frac{1}{2}+} } \\
& \leq ||\phi_3||_{X_{s, \frac{1}{2}+}}.
\end{align*}
It remains to bound $||\phi_4||_{L^3(\rr^{2+1})}$.  
Interpolating between
$||\phi_4||_{L^2_tL^2_x} \; \leq \; ||\phi_4||_{X_{0,0}}$
and the Strichartz estimate
$||\phi_4||_{L^4_tL^4_x} \; \lesssim \; ||\phi_4||_{\Xzohp}$
yields 
\begin{align*}
||\phi_4||_{L^3_tL^3_x} & \lesssim ||\phi_4||_{X_{0, \frac{1}{2} -- }}.
\end{align*}
This completes the proof of 
\eqref{dualwltdo}, and hence Proposition \ref{modifiedlocalexist}. 
\end{proof}

\begin{proof}[Proof of Proposition \ref{almostconservation}]
 The usual energy \eqref{energy} is shown to be conserved by
differentiating in time, integrating by parts, and using 
the equation \eqref{nls},
\begin{align*}
 \partial_t E(\phi) &= \Re \int_{\rr^2} \overline{\phi_t} (|\phi|^2 \phi -
\Delta \phi) dx \\
& = \Re \int_{\rr^2} \overline{\phi_t} ( |\phi|^2 \phi -
\Delta \phi - i \phi_t) dx\\
& = 0.
\end{align*} 
We follow the same strategy to estimate the growth of $E(I\phi)(t)$,
\begin{align*}
\partial_t E(I\phi)(t) & = \Re \int_{\rr^2} 
\overline{I(\phi)_t} ( |I\phi|^2 
I\phi - \Delta I\phi - iI\phi_t) dx \\
& = \Re \int_{\rr^2} \overline{I(\phi)_t} ( |I\phi|^2 I\phi -
I(|\phi|^2 \phi)) dx,
\end{align*}
where in the last step we've applied $I$ to \eqref{nls}. 
When
we integrate
in time and apply the Parseval formula\footnote{That is, 
$\int_{\rr^n} f_1(x) f_2(x) f_3(x) f_4(x) dx \; = \;
\int_{\xi_1 + \xi_2 + \xi_3 + \xi_4 = 0} \hat{f}_1(\xi_1) 
\hat{f}_2(\xi_2)\hat{f}_3(\xi_3)\hat{f}_4(\xi_4) $ where 
$\int_{\sum_i \xi_i = 0}$ here denotes integration with respect to the 
hyperplane's measure \newline 
$\delta_0(\xi_1+\xi_2+\xi_3+ \xi_4) d\xi_1 d\xi_2 d\xi_3 d\xi_4$, with
$\delta_0$
the one dimensional  Dirac mass.}
it remains for us to bound
\begin{multline}
E(I\phi(\delta))-E(I\phi(0))
  = \\
\int_0^\delta \int_{\sum_{j=1}^4 \xi_j = 0}
\left ( 1 - \frac{ m(\xi_2 + \xi_3 + \xi_4)}{m(\xi_2) 
\cdot m(\xi_3) \cdot m(\xi_4)} 
\right) \widehat{\overline{I \partial_t \phi}}
(\xi_1) \widehat{I\phi}(\xi_2) \widehat{\overline{I \phi}}(\xi_3)
\widehat{I \phi}(\xi_4). 
\label{fundII}
\end{multline}
The reader may ignore the appearance of complex conjugates here and in the sequel,
as they have no impact on the availability of estimates.  (See e.g. 
Lemma \ref{bourgainstrichartz} above.)  We include the complex conjugates 
for completeness.

We use the equation \eqref{nls} to substitute for 
$\partial_t I(\phi)$ in \eqref{fundII}.  Our aim is to show that 
\begin{align}
\termone + \termtwo & \lesssim N^{- \frac{3}{2}+},
\label{oneandtwo}
\end{align}
where the two terms on the left are 
\begin{align}
\termone & \equiv \left |\int_0^\delta \int_{\sum_{i=1}^4 \xi_i = 0}
\left( 1 - \frac{m(\xi_2 + \xi_3 + \xi_4)}{m(\xi_2)  m(\xi_3)
m(\xi_4)} \right) \widehat{( \Delta \overline{I\phi})}(\xi_1) 
\cdot \widehat{I\phi}(\xi_2)
\cdot \widehat{\overline{I\phi}}(\xi_3) \cdot \widehat{I\phi}(\xi_4) \right|
\label{I} \\
\termtwo & \equiv \left| \int_0^\delta \int_{\sum_{i=1}^4 \xi_i = 0}
\left( 1 - \frac{m(\xi_2 + \xi_3 + \xi_4)}{m(\xi_2)  m(\xi_3)
m(\xi_4)} \right) \widehat{(\overline{I(|\phi|^2 \phi))}}(\xi_1) 
\cdot \widehat{I\phi}(\xi_2)
\cdot \widehat{\overline{I\phi}}(\xi_3) \cdot \widehat{I\phi}(\xi_4) \right|.
\label{II} 
\end{align}
In both cases we break $\phi$ into a sum of dyadic constituents $\psi_j$, each
with frequency support $\langle \xi \rangle \thicksim 2^j$, 
$j =0, \ldots$.. 

For both $\termone$ and $\termtwo$ we'll pull the symbol 
\begin{equation}
\label{symbol}
1 -  \frac{m(\xi_2 + \xi_3 + \xi_4)}{m(\xi_2)  m(\xi_3)
m(\xi_4)} 
\end{equation} 
out of the integral, estimating it pointwise in absolute value, using 
two different strategies 
depending on the relative sizes of the frequencies
involved.  After so bounding the factor \eqref{symbol}, the remaining 
integrals in \eqref{I}, \eqref{II}, involving the pieces $\psi_i$ of 
$\phi$,
are estimated by reversing the Plancherel formula\footnote{
Assuming, as we may, that the spatial 
Fourier transform of $\phi$ is always
positive.} and 
using duality, H\"older's inequality, and Strichartz estimates. 
We can sum over the all frequency pieces $\psi_i$ since
our bounds decay geometrically in these frequencies.  We suggest that
the reader at first ignore this summation issue, and so ignore on first reading the appearance below of all
factors such as $N_i^{0-}$ which we include only to show explicitly why our frequency 
interaction estimates sum. The main goal of the analysis is to establish the decay of 
$N^{-\frac{3}{2} +}$ in each class of frequency interactions below.  

Consider first $\termone$. By Proposition \ref{modifiedlocalexist},
\begin{align*}
||\Delta (I \phi)||_{X^{\delta}_{-1, \frac{1}{2} +}} &
\leq ||I\phi||_{\Xdoohp} \\
& \lesssim 1.
\end{align*}
Hence we conclude $\termone \lesssim N^{-\frac{3}{2}+}$ once we show 
\begin{multline} \label{almostthere}
 \left|\int_0^\delta \int_{\sum_{i=1}^4 \xi_i =0} 
\left ( 1 - \frac{ m(\xi_2 + \xi_3 + \xi_4)}{m(\xi_2) \cdot m(\xi_3) 
\cdot m(\xi_4)} 
\right) \widehat{\overline{\phi_1}}(\xi_1) \widehat{\phi_2}(\xi_2) \widehat{ \overline{
\phi_3}}(\xi_3)
\widehat{\phi_4}(\xi_4) \right|  \\
\lesssim \quad N^{- \frac{3}{2}+} (N_1 N_2 N_3 N_4)^{0-}||\phi_1||_{X_{-1, \frac{1}{2}+}} 
\cdot || \phi_2||_{\Xoohp} \cdot ||\phi_3||_{\Xoohp} \cdot
||\phi_4||_{\Xoohp},
\end{multline}
for any functions $\phi_i, \, i = 1, \ldots, 4$ with positive spatial
Fourier transforms supported
on 
\begin{align}
\label{defineNi}
\langle \xi \rangle & \thicksim 2^{k_i} \quad \equiv N_i,
\end{align}
for some $k_i \in
\{0,1,\ldots\}$.  (Note that we are not decomposing the frequencies $|\xi| \leq 1$ here.
In the three dimensional argument we'll need to do this.)     
The inequality \eqref{almostthere} implies our desired bound 
\eqref{oneandtwo} for $\termone$ once we sum over
all dyadic pieces $\psi_j$.  
 
By the symmetry of the multiplier \eqref{symbol} 
in $\xi_2, \xi_3, \xi_4$, and the fact that the refined Strichartz
estimate \eqref{bstrich} allows complex conjugates on either factor, we may assume
for the remainder of this proof that 
\begin{equation}
\label{wemayassume}
N_2 \, \geq \, N_3 \, \geq \, N_4.
\end{equation}
Note too that $\sum_{i=1}^4 \xi_i = 0$ in the integration 
of \eqref{almostthere} so that $N_1 \lesssim N_2$.  Hence it is sufficient to obtain a 
decay factor
of $N^{- \frac{3}{2} +} N_2^{0-}$ on the right hand side of \eqref{almostthere}.
We now split the different frequency interactions into three cases,
according to the size of the parameter $N$ in comparison to the $N_i$.

\noindent{\bf $\termone$, Case 1:  $N \gg N_2$}.  According to
\eqref{Iproperties}, the symbol \eqref{symbol}
is in this
case identically zero and
the bound \eqref{almostthere} holds trivially.

\noindent{\bf $\termone$, Case 2:  $N_2 \gtrsim N \gg N_3 \geq N_4$}.  
Since $\sum_i \xi_i = 0$, we have here also $N_1 \thicksim N_2$.
By the mean value theorem,
\begin{align}
\left| \frac{m(\xi_2) - m(\xi_2 + \xi_3 + \xi_4)}{m(\xi_2)} \right|
& \lesssim
\frac{|\nabla m(\xi_2)\cdot(\xi_3 + \xi_4)|}{m(\xi_2)} 
\lesssim \frac{N_3}{N_2} .\label{pointwiseII}
\end{align}
This pointwise bound together with Plancherel's
theorem and \eqref{bstrich} yield
\begin{align}
\text{Left Side of \eqref{almostthere}} & \leq  \frac{N_3}{N_2}
||\phi_1 \phi_3 ||_{L^2([0, \delta] \times \rr^2])} 
||{\phi_2} \phi_4 ||_{L^2([0, \delta] \times \rr^2)} \\
&\leq \frac{N_3 N_3^{\frac{1}{2}} N_4^{\frac{1}{2}}}{N_2 N_1^{\frac{1}{2}}
N_2^{\frac{1}{2}}} 
\prod_i || \phi_i||_{\Xdzohp}.
\label{almostthereII}
\end{align}
Comparing \eqref{almostthere} with \eqref{almostthereII} it remains only 
to show that 
\begin{align*}
\frac{N_3 N_3^{\frac{1}{2}}N_4^{\frac{1}{2}}\langle N_1 \rangle}{N_2
N_1^{\frac{1}{2}}N_2^{\frac{1}{2}}N_2 \langle N_3\rangle \langle
N_4\rangle} &
\lesssim N^{- \frac{3}{2}+ } N_2^{0-}, 
\end{align*}
which follows immediately from our assumptions $N_1 \sim N_2
\gtrsim N \gg N_3 \geq N_4$.

\noindent{\bf $\termone$, Case 3: $N_2 \geq N_3 \gtrsim N$}. 
We use in this instance a trivial pointwise bound on the symbol,
\begin{align}
\label{trivial}
\left|1 - \frac{m(\xi_2 + \xi_3 + \xi_4)}{m(\xi_2) m(\xi_3) m(\xi_4)} \right|
& \lesssim \frac{m(\xi_1)}{m(\xi_2)m(\xi_3)m(\xi_4)}.
\end{align}
When estimating the remainder of the integrand on the left of
\eqref{almostthere}, break the interactions into
two subcases, depending on which frequency is comparable to $N_2$.

{\noindent {\bf Case 3(a):} $N_1 \sim N_2 \geq N_3 \gtrsim N$.} We aim 
for 
\begin{align*}
\frac{m(N_1)}{m(N_2) m(N_3) m(N_4)} \cdot  \left|
\int_{0}^{\delta}¥\int_{\sum_{i = 1}^4 \xi_i =0} \widehat{\overline{\phi_1}} 
\widehat{\phi_2} \widehat{\overline{\phi_3}} \widehat{\phi_4} \right|
& \lesssim \frac{N^{-\frac{3}{2}+} N_2^{0-}N_2 N_3 \langle N_4 \rangle}{N_1} 
\prod_{i=1}^4 
|| \phi_i ||_{\Xdzohp}.
\end{align*}
Pairing $\overline{\phi_1} \cdot \phi_4$ and $\phi_2 \cdot \overline{\phi_3}$ in 
$L^2$ and applying
\eqref{bstrich}, it remains to show 
\begin{align*}
\frac{m(N_1) N_4^{\frac{1}{2}} N_3^{\frac{1}{2}} }{m(N_2) m(N_3) m(N_4)
N_1^{\frac{1}{2}}N_2^{\frac{1}{2}}} & \lesssim N^{-\frac{3}{2}+}
\frac{ N_2^{1-} N_3 \langle N_4 \rangle}{N_1},
\end{align*}
or 
\begin{align}
\frac{N^{\frac{3}{2}+} N_2^{0-}}{m(N_3) m(N_4)N_2 N_3^{\frac{1}{2}}
 \langle N_4 \rangle^{\frac{1}{2}}} & \lesssim 1.
\label{getthis2d}
\end{align}
When estimating such fractions here and in the sequel, we frequently
use two trivial observations\footnote{Alternatively, use \eqref{Iproperties} to write out 
the value of $m$ explicitly.}: for any $p > \frac{3}{7}$, the function $m(x) x^p$ is 
increasing; 
and $m(x) \langle x \rangle ^p$ is bounded below.  For example, in the denominator of \eqref{getthis2d},
$m(N_4) \langle N_4 \rangle^{\frac{1}{2}} \gtrsim 1$ and $m(N_3)N_3^{\frac{1}{2}} \gtrsim 
m(N) N^{\frac{1}{2}} \, = \, N^{\frac{1}{2}}$.  After these observations
one quickly concludes that \eqref{getthis2d} holds.   

\noindent{{\bf Case 3(b):} $N_2 \sim N_3 \gtrsim N$}.  In this case we also 
know $N_1 \lesssim N_2$, since the frequencies $\xi_i$ must sum to zero.
Argue as above, now pairing $\phi_1 \phi_2$ and $\phi_3 \phi_4$ in $L^2$. 
The desired bound \eqref{almostthere} will follow from
\begin{align*}
\frac{m(N_1) N_1^{\frac{1}{2}} N_4^{\frac{1}{2}} }{m(N_2) m(N_3) m(N_4)
N_2^{\frac{1}{2}}N_3^{\frac{1}{2}}} & \lesssim N^{-\frac{3}{2}+}
\frac{ N_2^{1-} N_3 \langle N_4 \rangle}{\langle N_1 \rangle},
\end{align*}
or, after cancelling powers of $N_1$ in the numerator with 
powers of $N_2$ in the denominator,
\begin{align}
\frac{m(N_1) N^{\frac{3}{2}-} N_2^{0+}}{m(N_2) m(N_3) m(N_4) N_3^{\frac{1}{2}}
N_2 \langle N_4 \rangle^{\frac{1}{2}}} & \lesssim 1.
\label{getthisII2d}
\end{align}
Using $m(N_4)\langle N_4 \rangle^{\frac{1}{2}} \gtrsim 1$ and
that both $m(N_2) N_2^{\frac{1}{2}},\;  m(N_3) N_3^{\frac{1}{2}} 
\gtrsim m(N)N^{\frac{1}{2}} \; = \; N^{\frac{1}{2}}
$, 
we get \eqref{getthisII2d}.
This completes the proof of 
\eqref{almostthere}, and the bound for the contribution of 
$\termone$ in \eqref{oneandtwo}.

We turn to the bound \eqref{oneandtwo} for $\termtwo$
\eqref{II}.  As in our previous discussion of $\termone$, 
it suffices to show
\begin{multline}
\label{gottashowit}
\left| \int_0^{\delta} \int_{\sum_{i = 1}^6 \xi_i = 0} \left(1 - 
\frac{m(\xi_4 + \xi_5 + \xi_6)}{m(\xi_4) m(\xi_5) m(\xi_6)} \right)
P_{N_{123}}\widehat{\overline{I(\phi_1 \phi_2 \phi_3)}}(\xi_1+\xi_2+\xi_3)
\widehat{I\phi_4}(\xi_4) \widehat{\overline{I\phi_5}}(\xi_5)
\widehat{I\phi_6}(\xi_6) \right| \\
\; \lesssim \; \; N^{-\frac{3}{2}+}N_4^{0-} \prod_{i=1}^6 ||I\phi_i||_{\Xoohp},
\end{multline}
where as above, $0 \leq \widehat{\phi_i}(\xi_i)$ is supported for
$|\xi_i| \thicksim N_i = 2^{k_i}$, 
 and without loss of generality,
\begin{align}
\label{Nstuff}
N_4 \geq N_5 \geq N_6,\;  & \text{and} \;  N_4 \gtrsim N,
\end{align}
the latter assumption since otherwise the symbol on the 
left of \eqref{gottashowit} vanishes.
In \eqref{gottashowit} we have written $P_{N_{123}}$ for the projection onto 
functions supported in the $N_{123}$ dyadic spatial 
frequency shell. The decay factor on the right of \eqref{gottashowit}
allows us to sum in $N_4, N_5,N_6,$ and $N_{123}$, which suffices
as we do not dyadically decompose that part of $\termtwo$ represented
here by $\phi_i, i=1,2,3.$
We pointwise bound the symbol on the left of \eqref{gottashowit}  
in the obvious way 
\begin{align*}
\left| 1 - \frac{m(\xi_4 + \xi_5 + \xi_6)}{
m(\xi_4) m(\xi_5) m(\xi_6)} \right| 
& \lesssim\frac{m(N_{123})}{m(N_4)m(N_5)
m(N_6)}
\end{align*}
and as before, we undo the Plancherel formula.  After 
applying H\"older's inequality, it suffices to show
\begin{multline}
\label{gottashowitII}
\frac{m(N_{123})}{m(N_4) m(N_5) m(N_6)} \cdot 
||P_{N_{123}}I(\phi_1 \phi_2 \phi_3)||_{L^2_tL^2_x} \cdot
||I\phi_4||_{L^4_tL^4_x} \cdot || I \phi_5||_{L^4_tL^4_x}
\cdot ||I\phi_6||_{L^\infty_tL^\infty_x}  \\
 \lesssim \; \; N^{-\frac{3}{2}+}N_4^{0-} \prod_{i=1}^6 ||I\phi_i||_{\Xoohp}.
\end{multline}
To this end we'll use

\begin{lemma} \label{boundinstuff}
Suppose the functions $\phi_i, \; i =1,\ldots 6$ as above.
Then,
\begin{align}
\label{uno}
||P_{N_{123}} I(\phi_1 \phi_2\phi_3)||_{L^2_tL^2_x} & \lesssim
\frac{1}{\langle N_{123} \rangle} \prod_{i=1}^3 ||I\phi_i||_{\Xoohp}, \\
\label{dos}
||I\phi_j||_{L^4_tL^4_x} & \lesssim \frac{1}{\langle N_j \rangle} ||I\phi_j||_{
\Xoohp} \; \; j = 4,5,  \\
\label{tres}
||I\phi_6||_{L^\infty_t L^\infty_x} & \lesssim ||I\phi_6||_{\Xoohp}.
\end{align}
\end{lemma}
\begin{proof} 
For \eqref{uno}, it suffices to prove
\begin{align}
\label{unoshow}
||\langle \nabla  \rangle P_{N_{123}} I(\phi_1 \phi_2\phi_3)||_{L^2_tL^2_x} & \lesssim
\prod_{i=1}^3 ||I\phi_i||_{\Xoohp}.
\end{align}
(See Section 2 above for notation). The operator $\langle \nabla \rangle I$ obeys a Leibniz rule.  Using
H\"older's inequality on a typical resulting term,
\begin{align}
\label{firststep}
||P_{N_{123}} \left((\langle \nabla  \rangle I(\phi_1)) \phi_2\phi_3\right)||_{L^2_tL^2_x} & \lesssim
||\langle \nabla  \rangle I(\phi_1)||_{L^4_{x,t}}||\phi_2||_{L^8_{x,t}} ||\phi_3||_{L^8_{x,t}}.
\end{align}
By Sobolev's inequality and a $L^8_tL^{\frac{8}{3}}_x$ Strichartz
estimate \eqref{strichartz},
\begin{align}
||\phi_2||_{L^8_{t,x}} & \lesssim || \langle \nabla \rangle^{\frac{1}{2}} \phi_2 ||_{L^8_t
L^{\frac{8}{3}}_x} \nonumber \\
& \lesssim || \langle \nabla \rangle^{\frac{1}{2}} \phi_2 ||_{\Xzohp} \nonumber \\
& \lesssim || \phi_2 ||_{\Xoohp} \label{secondstep}
\end{align}
and similarly for the $\phi_3$ factor on the right of \eqref{firststep}.
Applying the $L^4_{x,t}$ Strichartz estimate,
\begin{align}
\label{thirdstep}
||\langle \nabla \rangle I\phi_1||_{L^4_{x,t}} & \lesssim ||I\phi_1||_{\Xoohp}.
\end{align}
Together, \eqref{firststep} - \eqref{thirdstep} yield \eqref{uno}.

The bounds \eqref{dos} follow immediately from the $L^4_{t,x}$ Strichartz
estimate as in \eqref{thirdstep}.  
The estimate \eqref{tres} is seen using Sobolev embedding,
the fact that $\phi_6$ is frequency localized, and the 
$L^{\infty}_t L^2_x$ Strichartz bound,
\begin{align*}
||I \phi_6||_{L^\infty_{x,t}} & \lesssim ||\langle \nabla \rangle I \phi_6||_{L^\infty_t
L^2_x} \\
& \lesssim ||I \phi_6||_{\Xoohp}.
\end{align*}
\end{proof}
Together, \eqref{gottashowitII} and Lemma \ref{boundinstuff}
leave us to show 
\begin{align}
\label{frodo}
\frac{m(N_{123}) \cdot N^{\frac{3}{2}-} N_4^{0+} }{m(N_4) m(N_5) m(N_6)
\langle N_{123} \rangle \langle N_4 \rangle \langle N_5 \rangle } & \lesssim 1
\end{align}
under the assumption \eqref{Nstuff}. We can
break the frequency interactions
into two cases:  $N_4 \thicksim N_5$ and $N_4 \thicksim N_{123}$, 
since we have $\sum_{i = 1}^6 \xi_i = 0 $ in \eqref{gottashowit}.

\noindent{\bf $\termtwo$, Case 1; $N_4 \thicksim N_5; N_4 \geq N_5 \geq
N_6; \; N_4 \gtrsim N$}:  We aim here for
\begin{align*}
\frac{m(N_{123})N^{\frac{3}{2}-} N_4^{0+}}{(m(N_4))^2 \langle N_4 \rangle^2 m(N_6) 
\langle N_{123} \rangle } &
\lesssim 1.
\end{align*}
Since 
$m(N_4)\langle N_4 \rangle ^{\frac{1}{2}} \, \gtrsim \, m(N) \langle N
\rangle^{\frac{1}{2}} \, =\, \langle N \rangle^{\frac{1}{2}}$
it suffices to show
\begin{align}
\frac{m(N_{123})N^{\frac{1}{2}-}N_4^{0+}}{  \langle N_4\rangle
  m(N_6) \langle N_{123}\rangle }
& \lesssim 1,
\end{align}
which is clear since $\langle N_{123}\rangle  \; \geq m(N_{123})$, and
\begin{align}
\label{increasing}
m(y)\langle x \rangle^{\frac{1}{2}} & \gtrsim 1 \quad \text{for all} \,  \quad 0 
\leq y \leq x. 
\end{align}

\noindent{\bf $\termtwo$, Case 2;  $N_4 \thicksim N_{123}; N_4 \geq 
N_5 \geq N_6; \; N_4 \gtrsim N$}:  Here we argue that
\begin{align*}
\frac{m(N_4) N^{\frac{3}{2}-} N_4^{0+} }{m(N_4) \langle N_4 \rangle^2 m(N_5) 
m(N_6) \langle N_5 \rangle}
&  \lesssim
\frac{N^{\frac{3}{2}-}N_4^{0+}}{m(N_5) \langle N_5\rangle^{\frac{1}{2}} m(N_6)
\langle N_4 \rangle^{2} \langle N_5\rangle^{\frac{1}{2}} } \lesssim 1,
\end{align*}
using \eqref{increasing} and our assumptions on the $N_i$.  
This completes the proof of \eqref{oneandtwo} and 
hence the proof of Proposition \ref{almostconservation}.
\end{proof}

\section{ Proof of Theorem 
\ref{maintheorem} in $\rr^3$}
\label{section:3d}
In three space dimensions our almost conservation law takes the following
form,

\begin{proposition} \label{almostconservation3d} 
Given $s > \frac{5}{6}, N \gg1, $ and initial data 
$\phi_0 \in C^{\infty}_0(\rr^3)$ 
with $E(I_N \phi_0) \leq 1$, then there exists a universal constant
$ \delta$ so that the solution
\begin{align*}
\phi(x,t) & \in C([0,\delta], H^s(\rr^3))
\end{align*}
of \eqref{nls}-\eqref{nlsdata} satisfies
\begin{align}
\label{increment3}
E(I_N\phi)(t) & =  E(I_N\phi)(0) + O(N^{- 1+}),
\end{align}
for all $t \in [0, \delta]$.
\end{proposition}

The norm $||\phi(t, \cdot)||_{L^2(\rr^3)}$ is supercritical
with respect to the scaling \eqref{scaling}.  Hence, aside from the $L^2$ conservation
\eqref{l2conservation} we will avoid using this quantity in the
proof of  the three dimensional result.  Beside the technical issues 
introduced by scaling the $L^2$ norm, our
proof of Theorem \ref{maintheorem} for $n=3$ follows 
very closely the $n=2$ arguments of Section \ref{section:2d}.  

We begin with the fact that Proposition \ref{almostconservation3d} implies
Theorem \ref{maintheorem} with $n=3$.  Recall that it suffices
to show the $H^s(\rr^3)$ norm of the solution to \eqref{nls}-
\eqref{nlsdata}
grows polynomially in time. Recall too $\philambda$ as the scaled 
solution defined in \eqref{scaling}.  When $n=3$, the definition 
of the energy \eqref{energy} and Sobolev embedding imply 
\begin{equation} \label{Ein3}
\begin{split}
E(I_N \philambda_0) & = \frac{1}{2} ||\nabla I_N \philambda_0 ||^2_{L^2(\rr^3)} 
 + \frac{1}{4} ||I_N \philambda_0 ||^4_{L^4(\rr^3)} \\
 & \leq C_0 N^{2 - 2s} \lambda ^{1 - 2s} 
 (1 + ||\phi_0||_{H^s(\rr^3)})^4.
\end{split}
\end{equation}
Once the parameter $N$ is chosen, we will choose $\lambda$
according to 
\begin{align}
\label{chooselambda3}
\lambda & = \left(\frac{1}{2C_0} \right)^{\frac{1}{1-2s}} N^{ \frac{2s-2}{1-2s}} (1 + ||\phi_0||_{H^s(\rr^3)} )^
{ - \frac{4}{1-2s}}.  
\end{align}
 Together, \eqref{chooselambda3} and \eqref{Ein3} give
 $E(I_N \philambda_0) \leq \frac{1}{2}$.
 We can therefore apply Proposition \ref{almostconservation3d} at least 
 $C_1\cdot N^{1-}$ times to give
\begin{align}
\label{Esize}
E(I_N \philambda)(C_1 N^{1-} \cdot \delta) & \thicksim 1.
\end{align}
The estimate \eqref{Esize} implies  $||\phi(t, \cdot)||_{H^s(\rr^3)}$
grows at most polynomially when $\frac{5}{6}<s<1$.  This can be seen
exactly as in the two dimensional case.  We include the argument here
for completeness.

Given any $T_0  \gg 1$, first choose $N \gg 1$ so that 
\begin{equation}
\label{chooseN}
T_0  = \frac{C_1 N^{1-} \delta}{\lambda^2} 
\thicksim N^{\left( \frac{ \frac{5}{2} - 3s - }
{\frac{1}{2} - s} \right) }.
\end{equation}
Note that the exponent of $N$ on the right of \eqref{chooseN} is positive
(and hence this definition of $N$ makes sense) precisely when $s > \frac{5}{6}$.  
In three space dimensions we have
\begin{align*}
E(I_N \philambda)(\lambda^2 t) & = \frac{1}{\lambda} E(I_N \phi)(t). 
\end{align*}
According to \eqref{chooselambda3}, \eqref{Esize}, \eqref{chooseN},
we therefore get
\begin{align*}
E(I_N \phi)(T_0) & \leq \lambda E(I_N \philambda)(\lambda^2 T_0) \\
& \lesssim \lambda \\
& \lesssim N^{\frac{2s-2}{1-2s}} \\
& \lesssim T_0^{\frac{1-s +}{3(s-\frac{5}{6})}}.
\end{align*}
According to \eqref{boundhsnorm} and \eqref{l2conservation}, 
the $H^s(\rr^3)$ norm grows
with at most half this rate when $\frac{5}{6} < s < 1$,
\begin{align}
\label{polynomial3d}
||\phi||_{H^s(\rr^3)}(T) & \lesssim (1 +
T)^{\frac{1-s +}{6(s-\frac{5}{6})}}.
\end{align}

As in the two dimensional argument, the proof
of Proposition \ref{almostconservation3d} relies on bounds 
for the local-in-time $H^s$ solution. The following analogue
of Proposition \ref{modifiedlocalexist} avoids the use of 
the norm $||\phi(\cdot, t)||_{L^2(\rr^3)}$, which, as mentioned
above, is supercritical with respect to scaling.

\begin{proposition} \label{modifiedlocalexist3d} Assume $\frac{5}{6} <
s < 1$ and we are given data for \eqref{nls}-\eqref{nlsdata} with
$E(I\phi_0) \; \leq \; 1$.  Then 
there is a universal constant $\delta >0$
so that the solution $\phi$ obeys
the following bound on the time interval $[0, \delta]$,
\begin{align}
\label{localbounds3d}  
||\nabla I\phi||_{\Xdzohp} 
& \lesssim 1.
\end{align}
\end{proposition}

\begin{proof}
Arguing as in the proof of Proposition \ref{modifiedlocalexist},
it suffices to prove
\begin{align*}
||\nabla I(\phi \bar{\phi}\phi) ||_{X^{\delta}_{0, - \frac{1}{2} ++}}
& \lesssim ||\nabla I \phi ||^3_{\Xzohp},
\end{align*}
Again, the interpolation lemma in \cite{iteamIII}
allows us to assume
$N = 1$ in the definition \eqref{Ioperator} of the operator $I$.
After applying a Leibniz rule for the operator $\nabla I$ and
duality, we aim to show
\begin{align}
\label{goal3d}
||(\nabla I)(\phi_1) \cdot \overline{\phi_2} \cdot \phi_3 \cdot \psi
||_{L^1(\rr^{3+1})} & \lesssim ||\psi||_{X_{0, \frac{1}{2}--}}
\prod_{i=1}^3 ||\nabla I \phi_i||_{\Xzohp}.
\end{align}
Again, the complex conjugate will have no bearing on our bounds.  
We split the functions $\phi_j, \, j=2,3$ into high and low frequency 
components,
\begin{align}
\label{highlow}
\phi_j& = \phihigh_j + \philow_j,
\end{align}
where 
\begin{align*}
\text{supp} \, \widehat{\phihigh_j}(\xi, t) & \subset
 \{|\xi| \geq \frac{1}{2} \} \\
 \text{supp} \, \widehat{\philow_j}(\xi, t) & \subset
 \{|\xi| \leq 1 \}. 
 \end{align*}
 Note that when $n=3$, homogeneous Sobolev embedding
 and the $L^{10}_t L^{\frac{30}{13}}_x$ Strichartz estimate
 give 
 \begin{equation}
 \label{L10estimate}
 \begin{split}
 ||\phi||_{L^{10}_tL^{10}_x(\rr^{3+1})} & \lesssim 
 ||\nabla \phi||_{L^{10}_tL^{\frac{30}{13}}_x(\rr^{3+1})} \\
 & \lesssim ||\nabla \phi ||_{\Xzohp}.
 \end{split}
 \end{equation} 
 Consider first the low frequency components on the left
 of \eqref{goal3d}.  Apply Holder's inequality with the
 factors in $L^{\frac{10}{3}}_{x,t}, L^{10}_{x,t}, L^{10}_{x,t},$ and
 $L^2_{x,t}$ respectively.  The $L^{\frac{10}{3}}_{x,t}$ Strichartz
 estimate along with \eqref{L10estimate} give,
 \begin{align*}
 ||(\nabla I)(\phi_1) \cdot \overline{\philow_2} \cdot \philow_3 \cdot \psi
||_{L^1(\rr^{3+1})} & \lesssim ||\psi||_{X_{0,0}} ||\nabla I \phi_1||_{
X_{0, \frac{1}{2}+}} \prod_{i=2}^3 ||\nabla \philow_j||_{\Xzohp}.
\end{align*}
Together with the fact that $\philow_j = I\philow_j$, this bound accounts
for part of the low frequency contributions of $\phi_2, \phi_3$
in \eqref{goal3d}.  A typical contribution which remains 
to be bounded is  
$$
||\nabla(I\phi_1) \overline{\phihigh_2} \phihigh_3 \psi||_{L^1_{x,t}(\rr^{3+1})}.
$$
Recall the Strichartz estimate,
\begin{align}
||\psi||_{L^{\frac{10}{3}}_{x,t}(\rr^{3+1})} & \lesssim 
||\psi||_{\Xzohp}.  \label{L10/3estimate}
\end{align}
Interpolating between \eqref{L10/3estimate} and the trivial bound  
$||\psi||_{L^2_{x,t}(\rr^{3+1})} \lesssim ||\psi||_{X_{0,0}}$ gives
\begin{align}
||\psi||_{L^{3}_{x,t}(\rr^{3+1})} & \lesssim ||\psi||_{X_{0,
\frac{1}{2}--}}.
\label{psiL3}
\end{align}
Using \eqref{psiL3} and Holder's inequality on the left of 
\eqref{goal3d}, we aim to show
\begin{align}
\label{remains3d}
||\nabla(I\phi_1) \phihigh_2 \phihigh_3||_{L^{\frac{3}{2}}
_{x,t}(\rr^{3+1})} & \lesssim 
\prod_{i=1}^3 ||\nabla I \phi_i||_{\Xzohp}.
\end{align}
Since we've reduced to the case $N=1$, we note $I^{-1} = \langle \nabla 
\rangle^g$,
where  
\begin{align}
\label{thegap}
g & \equiv 1 - s \quad \in \; (0, \frac{1}{6})
\end{align}
is the gap between $s$ and $1$.    
We may therefore rewrite our desired estimate as 
 \begin{align}
\label{remains3dII}
||\nabla(I\phi_1) (\fatD^g I \overline{\phihigh_2}) (\fatD^g 
I \phihigh_3)||_{L^{\frac{3}{2}}
_{x,t}(\rr^{3+1})} & \leq 
\prod_{i=1}^3 ||\nabla I \phi_i||_{\Xzohp}.
\end{align}
But this estimate follows after taking the factors on the left
in $L^{\frac{10}{3}}_{x,t}, \, L^{\frac{60}{11}}_{x,t}, \, L^{
\frac{60}{11}}_{x,t}$, respectively and using H\"olders inequality.  
The first resulting factor is bounded using the $L^{\frac{10}{3}}_{x,t}$
Strichartz estimate.  As for the second two factors,  
Sobolev embedding, the bounds \eqref{thegap} on $g$, and
the 
$L^{\frac{60}{11} \, \frac{45}{17}}_{t,x}$ Strichartz estimate yield
for $j = 2,3$,
\begin{align*}
||\fatD^g I \phihigh_j||_{L^{\frac{60}{11}}_{x,t}(\rr^{3+1})} & 
\lesssim || \fatD^{1-g} \fatD^g I \phihigh_j||_
{L^{\frac{60}{11} \, \frac{45}{17}}_{t,x}} \\
& \lesssim ||\nabla I \phihigh_j ||_{\Xzohp}.
\end{align*}
The case where $\philow_2 \phihigh_3$ appears on the left of 
\eqref{goal3d} is handled similarly, using a homogeneous Sobolev
embedding to bound the $\philow_2$ term.  
\end{proof} 
\begin{proof}[Proof of Proposition \ref{almostconservation3d}]
Arguing as in the two dimensional result leaves us to show
\begin{align}
\label{oneandtwo3d}
\termone + \termtwo & \lesssim N^{- 1++},
\end{align}
where the two terms on the left are as before, \eqref{I}, \eqref{II}.
We will have to pay closer attention here than in $\rr^2$ when we sum the 
various dyadic components of this
estimate.  The fact that we only control inhomogeneous
norms \eqref{localbounds3d} forces us to decompose the frequencies
$|\xi| \leq 1$  as well.  

Considering first $\termone$, it follows from the definition
of the $X_{s,b}$ norms \eqref{norms} that 
\begin{align}
||\Delta I \phi||_{X_{-1, \frac{1}{2}+}^\delta} & \lesssim
|| \nabla I \phi ||_{X_{0, \frac{1}{2}+}^\delta}.
\end{align}
We conclude $\termone \lesssim N^{-1++}$ once we prove
\begin{multline} \label{almostthere3d}
 \left| \int_0^\delta \int_{\sum_{i=1}^4 \xi_i = 0}
\left ( 1 - \frac{ m(\xi_2 + \xi_3 + \xi_4)}{m(\xi_2) \cdot m(\xi_3) \cdot m(\xi_4)} 
\right) \widehat{\phi_1}(\xi_1) \widehat{\phi_2}(\xi_2) \widehat{ \phi_3}(\xi_3)
\widehat{\phi_4}(\xi_4) \right|  \\
\lesssim \quad N^{- 1++} C(N_1, N_2, N_3, N_4) ||\phi_1||_{X_{-1, \frac{1}{2}+}} 
\cdot \prod_{j=2}^4|| \nabla \phi_j||_{\Xzohp}, 
\end{multline}
for sufficiently small $C(N_1, N_2, N_3, N_4)$ and for any smooth 
functions $\phi_i, \, i = 1, \ldots 4$ with
$0 \leq \widehat{\phi_i}(\xi_i)$ supported
for $|\xi_i| \, \thicksim \, N_i \equiv 2^{k_i}, \, k_i = 0, \pm1, \pm2, \ldots$.  As before, we 
may assume $N_2 \geq N_3 \geq N_4$.  The precise extent to which 
$C(N_1, N_2, N_3, N_4)$ decays in its arguments, and the fact
that this decay allows us to sum over all dyadic shells, will be described
below on a case-by-case basis.  

In addition to the estimates \eqref{L10estimate}, \eqref{L10/3estimate}, our
analysis here uses the following related bounds, all of which are quick consequences 
of homogeneous Sobolev embedding, H\"older's inequality in the time variable, and/or 
Strichartz estimates. These estimates will allow for bounds decaying in the frequencies.
For a function $\phi$ with frequency
support in the $D$'th dyadic shell,
\begin{align}
||\phi||_{L^{10}_tL^{10\pm}_x ([0, \delta] \times \rr^3)} & 
\lesssim D^{0\pm} || \nabla \phi ||_{\Xzohp} \label{L10estimateplusminus} \\
||\phi||_{L^{\frac{10}{3}}_tL^{\frac{10}{3} -}_x([0, \delta] \times \rr^3)} & 
\lesssim \delta^{0+}  || \phi ||_{\Xzohp} \label{L10/3estimateminus} \\
||\phi||_{L^{\frac{10}{3}}_tL^{\frac{10}{3} +}_x([0, \delta] \times \rr^3)} & 
\lesssim D^{0+} || \phi ||_{\Xzohp}. \label{L10/3estimateplus}
\end{align}
\noindent{\bf $\termone$, Case 1:  $N \gg  N_2$}.  Again, 
the symbol \eqref{symbol}
is in this
case identically zero and
the bound \eqref{almostthere3d}  holds trivially, with $C \equiv 0$. \newline  
\noindent{\bf $\termone$, Case 2:  $N_2 \gtrsim N \gg  N_3 \geq N_4$}.
We have $N_2 \thicksim N_1$ here as well.  
We will show 
\begin{align}
\label{Ccase2}
C(N_1, N_2, N_3, N_4) & =
N_2^{0-} N_4^{0+}.
\end{align} 
With this decay factor, and the fact that we are considering here terms
where $N_1 \thicksim N_2$, we may immediately sum over the $N_1, N_2$ indices.
Similarly, the factor $N_4^{0+}$ in \eqref{Ccase2} allows us to sum over all 
terms here with
$N_3, N_4 \ll 1$.  It remains to sum the terms where 
$1 \lesssim N_4 \leq N_3 \ll N$, but these introduce
at worst a divergence  $N^{0+}\log(N)$, which is absorbed by the decay factor
$N^{-1++}$ on the right side of \eqref{almostthere3d}.  

We now show \eqref{almostthere3d}, \eqref{Ccase2}.  As before, \eqref{pointwiseII},
we bound the symbol in this case by $\frac{N_3}{N_2}$.  
We apply H\"older's inequality to the left side of \eqref{almostthere3d},
bounding $\phi_1, \phi_3$ in $L^{\frac{10}{3}}_{x,t}$ as in \eqref{L10/3estimate};
$\phi_2$ in $L^{\frac{10}{3}}_tL^{\frac{10}{3}-}_x$ as in \eqref{L10/3estimateminus}; and 
$\phi_4$ in $L^{10}_tL^{10+}_x$
as in \eqref{L10estimateplusminus} to get, 
\begin{align*}
\text{Left Side of \eqref{almostthere3d}} & \lesssim 
N_4^{0+}\frac{N_3}{N_2} ||\phi_1||_{\Xzohp} ||\phi_2||_{\Xzohp} 
||\phi_3||_{\Xzohp} ||\nabla \phi_4||_{\Xzohp} \nonumber \\
& \lesssim \frac{N_4^{0+} N_3  N_1 }{N_2 \cdot N_2 \cdot N_3}
||\phi_1||_{X_{-1, \frac{1}{2}+}} \prod_{j=2}^4 ||\nabla \phi_j||_{
\Xzohp}. 
\end{align*}
We conclude the bound \eqref{almostthere3d}, \eqref{Ccase2} for this case once we note
\begin{align*}
\frac{N_3 N_1 N^{1-} N_2^{0+}}{N_2 N_2 N_3} & \lesssim 1,
\end{align*}
which is immediate from our assumptions on the $N_i$. \newline
\noindent{\bf $\termone$, Case 3, $N_2 \geq N_3 \gtrsim N$:}
As in the two dimensional argument,  we use here the straightforward 
bound \eqref{trivial}
for the symbol.  The estimate of the remainder of the integrand
will break up into six different subcases, depending on which $N_i$ is
comparable to $N_2$, and whether or not $N_1, N_4 \ll  1$. \newline 
{\noindent {\bf Case 3(a)}, $N_1 \thicksim N_2 \geq N_3 \gtrsim N; N_4 \ll 1$:}
We will show here
\begin{align}
\label{Ccase3a}
C(N_1, N_2, N_3, N_4) & =
N_4^{0+} N_3^{0-},
\end{align} 
which suffices since one may use \eqref{Ccase3a} to sum directly in $N_3, N_4$, and use 
Cauchy-Schwarz to sum
in $N_1, N_2$.   

To establish \eqref{almostthere3d}, \eqref{Ccase3a}, estimate the $\phi_4$ factor 
in $L^{10}_tL^{10+}_x$ using \eqref{L10estimateplusminus};
$\phi_3$ in $L^{\frac{10}{3}}_t L^{\frac{10}{3}-}_x$ as in \eqref{L10/3estimateminus};
and $\phi_1, \phi_2$ in $L^{\frac{10}{3}}_{x,t}$ as in \eqref{L10/3estimate}.
It remains then to show
\begin{align}
\frac{m(N_1)N_1 N^{1-} N_3^{0+}}{m(N_2) m(N_3) N_2 N_3} & \lesssim 1 . \label{case3agoal}
\end{align}
Note that since $s \in (\frac{5}{6}, 1)$, we can use the following
fact while working in three space dimensions,
\begin{equation}
\label{3drule}
(m(x))^{p_1} x^{p_2} \quad \text{is nondecreasing in $x$ when} \,
0\, < \, p_1 \, \leq \, 6 p_2.
\end{equation}
We check \eqref{case3agoal} by first cancelling factors
involving $N_1$ and $N_2$ from numerator and denominator,and then
using \eqref{3drule},
\begin{align*}
\frac{m(N_1)  N_1  N^{1-} N_3^{0+}}{m(N_2)m(N_3) N_2 N_3}
& \lesssim \frac{N^{1-} N_3^{0+}}{m(N_3) N_3} \\
& \lesssim \frac{N^{1-} N_3^{0+}}{(m(N)) N^{1-} N_3^{0+}} \quad \lesssim 1. 
\end{align*}
%
%
%
{\noindent {\bf Case 3(b)}, $N_2 \sim N_3 \gtrsim N, N_1 \gtrsim 1, N_4 \ll 1$:}
Exactly as above, one shows \eqref{almostthere3d}, \eqref{Ccase3a} holds.
With this, one may sum directly in $N_4$, and also in  $N_1, N_2, N_3$ using the $N_3$ decay in 
\eqref{Ccase3a}.
\newline {\noindent {\bf Case 3(c)}, $N_2 \sim N_3 \gtrsim N, N_1 \ll 1, N_4 \ll 1$:}
Here we have 
\begin{align}
\label{Ccase3c}
C(N_1, N_2, N_3, N_4) & =
N_3^{0-} N_4^{0+}N_1^{0+}.
\end{align}
Allowing us to sum directly in all $N_i$.   One shows \eqref{Ccase3c}  
by modifying the argument in 3(a), taking 
$\phi_1, \phi_2$ in 
$L^{\frac{10}{3}}_t L^{\frac{10}{3}+}_x, L^{\frac{10}{3}}_t L^{\frac{10}{3}+}_x$,
respectively. 
\newline{\noindent {\bf Case 3(d)}, $N_1 \sim N_2 \geq N_3 \gtrsim N; N_4 \gtrsim 1$:}
We will show here
\begin{align}
\label{Ccase3d}
C(N_1, N_2, N_3, N_4) & =
N_3^{0-} N_4^{0-},
\end{align} 
allowing us to sum immediately in $N_3$ and $N_4$; summing in $N_1, N_2$ using Cauchy-Schwarz.

After taking the symbol out of the left side of \eqref{almostthere3d}
using \eqref{trivial}, we apply H\"older's inequality as
follows: estimate the $\phi_4$ factor in $L^{10}_tL^{10-}_x$ using \eqref{L10estimateplusminus};
$\phi_3$ in $L^{\frac{10}{3}}_t L^{\frac{10}{3}+}_x$ as in \eqref{L10/3estimateplus};
and $\phi_1, \phi_2$ in $L^{\frac{10}{3}}_{x,t}$ as in \eqref{L10/3estimate}.
We will establish \eqref{almostthere3d}, \eqref{Ccase3d} once we show
\begin{align}
\frac{m(N_1)N_1 N^{1--}N_3^{0+} N_3^{0+}}{m(N_2) m(N_3) m(N_4) N_2 N_3} & \lesssim 1 . \label{case3bgoal}
\end{align} 
This is done as in the argument of Case 3(a).  
\newline{\noindent {\bf Case 3(e)}, $N_2 \sim N_3 \gtrsim N$; $N_4 \gtrsim 1$;$N_1 \gtrsim 1$:}
We will show here
\begin{align}
\label{Ccase3e}
C(N_1, N_2, N_3, N_4) & =
N_2^{0-} N_1^{0-},
\end{align} 
allowing us to sum directly in all the $N_i$.
The H\"older's inequality argument here takes 
the $\phi_1$ factor in $L^{10}_tL^{10-}_x$ using \eqref{L10estimateplusminus};
$\phi_2$ in $L^{\frac{10}{3}}_t L^{\frac{10}{3}+}_x$ as in \eqref{L10/3estimateplus};
and $\phi_3, \phi_4$ in $L^{\frac{10}{3}}_{x,t}$ as in \eqref{L10/3estimate}.
We will have shown \eqref{almostthere3d}, \eqref{Ccase3e} once we show,
\begin{align}
\frac{m(N_1)N_1 N^{1--}N_2^{0+} N_2^{0+}}{m(N_2) m(N_3) m(N_4) N_2 N_3} & \lesssim 1 . \label{case3egoal}
\end{align} 
This argument is by now straightforward,
\begin{align*}
\frac{m(N_1)N_1 N^{1--}N_2^{0+} N_2^{0+}}{m(N_2) m(N_3) m(N_4) N_2 N_3}
& \lesssim \frac{N^{1--}N_2^{0+} N_2^{0+}}{m(N_2)^2  N_2} \\
& \lesssim \frac{N^{1--}N_2^{0+} N_2^{0+}}{ N^{1--} N_2^{0+} N_2^{0+}} \quad \lesssim 1.
\end{align*}
{\noindent {\bf Case 3(f)}, $N_2 \sim N_3 \gtrsim N$; $N_4 \gtrsim 1$;$N_1 \ll 1$:}
We will show here
\begin{align}
\label{Ccase3f}
C(N_1, N_2, N_3, N_4) & =
N_2^{0-} N_1^{0+},
\end{align} 
allowing us to sum directly in all the $N_i$.

The proof of \eqref{almostthere3d}, \eqref{Ccase3f}  
is similar to Case (3e), now taking $\phi_1$ in 
$L^{10}_tL^{10+}_x$ using \eqref{L10estimateplusminus};
$\phi_2$ in $L^{\frac{10}{3}}_t L^{\frac{10}{3}+}_x$ as in \eqref{L10/3estimateplus};
and $\phi_3, \phi_4$ in $L^{\frac{10}{3}}_{t}L^{\frac{10}{3}}_x$ as in \eqref{L10/3estimateminus}.

This completes the 3-dimensional analysis of $\termone$ in 
\eqref{oneandtwo3d}.

We will show $\termtwo \, \lesssim \, N^{-1+}$ using the straightforward
bound \eqref{trivial} on the symbol in the case $N_2 \geq N$, 
and the following,

\begin{lemma}  \label{lemma:cubebound}
\begin{align}
\label{cubebound}
||I(\phi_1 \phi_2 \phi_3)||_{L^2_tL^2_x([0,\delta] \times \rr^{3})}
& \lesssim \prod_{i=1}^3 ||\nabla I\phi_i||_{\Xzohp}.
\end{align}
\end{lemma}

We postpone the proof of Lemma \ref{lemma:cubebound}.  As in the work for $\termone$ above,
the argument bounding $\termtwo$ is complicated only by the presence of low frequencies.
Our aim is to show
\begin{align} \label{gottashowitIII}
\text{Left Side of \eqref{gottashowit}} & \lesssim C(N_{123}, N_4, N_5, N_6) 
\prod_{i=1}^6 ||\nabla I\phi_i||_{\Xzohp}. 
\end{align}
where $N_4 \geq N_5 \geq N_6$ and $N_4 \gtrsim N$, and as in the $\termone$ work
above, $C(N_{123},N_4,N_5,N_6)$ decays sufficiently fast to allow us to add up the
individual frequency interaction estimates to get \eqref{oneandtwo3d}.  

We sum first the interactions involving  $N_i \gtrsim 1$ for all frequencies 
in \eqref{gottashowitIII}.  In this case we'll show a decay factor of $C = N^{-1+} (N_{123}N_4N_5N_6)^{0-}
$, allowing us to sum in each index $N_i$ directly.  Apply H\"older's
inequality to 
the integrand on the left of \eqref{gottashowit}, taking the factors
in $L^2_{x,t}; L^{\frac{10}{3}}_{x,t};
L^{10}_{x,t}$; and $L^{10}_{x,t}$; respectively.  Using \eqref{trivial}, \eqref{cubebound},
the Strichartz estimates and \eqref{L10estimate} as 
in the $\termone$ argument, it suffices to show
\begin{align}
\label{laststop}
\frac{N^{1-} N_4^{0+} m(N_{123})}{N_4 m(N_4)m(N_5)m(N_6)} & \lesssim 1.
\end{align}
The fact that $m(x)$ is nonincreasing in $x$ and \eqref{3drule}
give us
\begin{align*}
\frac{N^{1-} N_4^{0+} m(N_{123})}{N_4 m(N_4)m(N_5)m(N_6)} & \lesssim
\frac{N^{1-}N_4^{0+}  m(N_{123})}{N_4 (m(N_4))^3}  \\
& \lesssim \frac{N^{1-} N_4^{0+} m(N_{123})}{N_4^{0+} N^{1-} (m(N))^3} \quad \lesssim 1.
\end{align*}
The above argument is easily modified in the presence of small frequencies.
We sketch these modifications here.  In case $N_{123} \sim N_4$, with $N_6 \ll 1$ and 
possibly also $N_5 \ll 1$, we need to get factors of $N_6^{0+}$ and possibly also
$N_5^{0+}$  on the right hand side of \eqref{gottashowitIII}.  We accomplish this by 
taking the factor $I\phi_6$ and possibly also $I\phi_5$ in  
$L^{10}_{t}L^{10+}_x$,and take the factor $I\phi_4$ in 
$L^{\frac{10}{3}}_t L^{\frac{10}{3}-}_x$, or possibly $L^{\frac{10}{3}}_t L^{\frac{10}{3}--}_x$.

In case $N_4 \sim  N_5$, with $N_{123}$ and/or $N_6$ small, a similar argument
gets the necessary decay:  we can take $P_{N_{123}}I(\phi_1 \phi_2 \phi_3)$ in $L^{2}_tL^{2+}_x$, and/or 
$I\phi_6$
in $L^{10}_tL^{10+}_x$, and take $I\phi_4$ in
$L^{10/3}_tL^{\frac{10}{3}-}_x$ or $L^{10/3}_tL^{\frac{10}{3}--}_x$.
\end{proof}    
\begin{proof}[Proof of Lemma \ref{lemma:cubebound}]
By the interpolation lemma  in
\cite{iteamIII}, we may assume $N=1$. By Plancherel's theorem,
it suffices to prove
\begin{align}
\label{sufficescubic}
|| \phi_1 \cdot \phi_2 \cdot \phi_3 ||_{L^2_{x,t}
([0,\delta]\times \rr^3)} & 
\lesssim \prod_{i=1}^3 ||\nabla \phi_i||_{\Xzohp}.
\end{align}  
Decomposing $\phi_i = \philow_i + \phihigh_i$ as in \eqref{highlow},
we consider first the contribution when only the low frequencies interact
with one another.  H\"older's inequality in space-time, 
homogeneous Sobolev embedding, H\"older's inequality in time,
and the energy estimate yield,
\begin{align*}
||\philow_1 \cdot \philow_2 \cdot 
\philow_3||_{ L^2_{x,t}([0,\delta]\times \rr^3)} & = 
||I \philow_1 \cdot I \philow_2 \cdot 
I\philow_3||_{L^2_{x,t}([0,\delta]\times \rr^3)} \\
& \lesssim \prod_{i=1}^3||I \philow_i||_{L^6_{x,t}([0,\delta]\times \rr^3)} \\
& \lesssim \prod_{i=1}^3||\nabla I \philow_i||_{L^6_tL^2_x
([0,\delta]\times \rr^3)}\\
& \lesssim \delta^{\frac{1}{2}} \prod_{i=1}^3 ||\nabla I \philow_i||_
{L^\infty_tL^2_x([0,\delta]\times \rr^3)} \\
& \lesssim  \prod_{i=1}^3 ||\nabla I \philow_i||_{\Xzohp}.
\end{align*}
A typical term whose contribution to \eqref{sufficescubic} remains
to be controlled is 
\begin{equation}
\label{remainscubic3d}
||I \philow_1 \cdot \fatD^g I \phihigh_2 \cdot \fatD^g I \phihigh_3 
||_{L^2_{x,t}([0,\delta]\times \rr^3)},
\end{equation}
where $g$ is as in \eqref{thegap}.  Take the first factor here
in $L^{6}_t L^{18}_x$ and each of the second two in 
$L^{6}_tL^{\frac{9}{2}}_x$ via H\"older's inequality.  Note
then that Sobolev embedding and the $L^6_tL^{\frac{18}{7}}_x$
Strichartz inequality give us
\begin{align*}
||I \philow_1||_{L^6_tL^{18}_x(\rr^{3+1})} & \lesssim
||\nabla I\philow_1 ||_{L^6_tL^{\frac{18}{7}}_x(\rr^{3+1})} \\
& \lesssim ||\nabla I \philow_1 ||_{\Xdzohp}.
\end{align*}
Similarly, the fact that $g \in (0, \frac{1}{6})$, Sobolev
embedding, H\"older's inequality in time, and the energy estimate
give us for $j=2,3$,
\begin{align*}
||\fatD^g I \phihigh_j||_{L^6_tL^{\frac{9}{2}}_x([0,\delta]\times \rr^3)}
& \lesssim ||\fatD I \phihigh_j ||_{L^6_t L^2_x([0,\delta]\times \rr^3)} \\
& \lesssim \delta^{\frac{1}{6}} 
||\nabla I \phihigh_j ||_{L^\infty_t L^2_x(\rr^{3+1})} \\
& \lesssim ||\nabla I \phi_j||_{\Xdzohp}.
\end{align*}
This completes the proof of Lemma \ref{lemma:cubebound}.
\end{proof}

\end{document}